\documentstyle[draft,amscd,syntonly,amssymb]{amsart}

\theoremstyle{plain}
\newtheorem{Thm}{Theorem}
\newtheorem{Prop}[Thm]{Proposition}

 \theoremstyle{definition}
\theoremstyle{remark}

\numberwithin{equation}{section}

\begin{document} 
 \title{On the homotopy of symplectomorphism groups of homogeneous spaces}

 \author{ ANDR\'{E}S   VI\~{N}A}
\address{Departamento de F\'{i}sica. Universidad de Oviedo.   Avda Calvo
 Sotelo.     33007 Oviedo. Spain. } 
 \email{avinae@@correo.uniovi.es}

 \maketitle
\begin{abstract}
Let ${\cal O}$ be a quantizable coadjoint orbit of a semisimple Lie group $G$.
Under certain hypotheses we prove that 
$\#(\pi_1(\text{Ham}({\cal O})))\geq \#(Z(G))$, where $\text{Ham}({\cal O})$
is the group of Hamiltonian symplectomorphisms of ${\cal O}$.

\end{abstract}
\section {Introduction} \label{S:intro}

 Let $(M,\omega)$ be a quantizable symplectic manifold \cite{nW92}. By $\text{Ham}(M,\omega)$
we denote the group of Hamiltonian symplectomorphisms of $(M,\omega)$ \cite{Mc-S}.
The group $\pi_1(\text{Ham}(M,\omega))$ is completely known in very rare cases
(see \cite[page 52]{lP01} \cite{McD}). In this note we determine a lower bound
for the cardinal of this homotopy group, when $M$ is an homogeneous space.
Our approach is based in  properties of the representation ``symplectic action".
In   \cite{aV02} we proved the existence a representation $\kappa$ of
the group $\pi_1(\text{Ham}(M,\omega))$. $\kappa$ associates to $[\psi]$ the 
symplectic action around the loop $\psi$ in $\text{Ham}(M,\omega)$. We also established 
some properties of $\kappa([\psi])$, which permit us to calculate it in particular cases. 
If $G$ is a semisimple Lie group and the quantizable manifold is the
 coadjoint orbit  ${\cal O}_{\eta}$ of an 
element $\eta\in{\frak g}^*$, each curve $\{g_t\in G\}_{t\in[0,1]}$ in $G$ with $g_0=e$
and $g_1\in Z(G)$ determines a loop $\psi$ in $\text{Ham}({\cal O}_{\eta})$.  In this case
$\kappa([\psi])$ is the value at $g_1$ of a character of the stabilizer $G_{\eta}$ of $\eta$
for the coadjoint action of $G$. Under certain hypotheses, we prove that  such two curves
 $\{g_t\}$ and $\{\tilde g_t\}$ with  different endpoints generate two loops $\psi$ and
$\tilde\psi$, such that $\kappa([\psi])\ne \kappa([\tilde\psi])$. Hence
$\#(\pi_1(\text{Ham}({\cal O}_{\eta})))\geq\#(Z(G))$.

In Section 2 the definition of the representation $\kappa$ and its general properties
 are reviewed.
In Section 3 we determine a lower bound to $\#(\pi_1(\text{Ham}(M,\omega)))$ when
$(M,\omega)$ is a coadjoint orbit. As particular case
coadjoint orbits of $SU(n)$ are considered, and we prove that
$\pi_1(\text{Ham}({\Bbb C}P^{n-1},\,\omega_{\text{Fubini-Study}}))$ has at least
$n$ elements.

%%%%%%%%%%%%%%%%%%%%%%%%%%%%%%%%%%%%%%%%%%%%%%%%%%%%%%%%%%%%%%%%%%%%%%%%%%%%%
%%%%%%%%%%%%%%%%%%%%%%%%%%%%%%%%%%%%%%%%%%%%%%%%%%%%%%%%%%%%%%%%%%%%%%%%%%%%%%%%

  \smallskip

\section{The representation $\kappa$}\label{S:Repre} 

Let $(M,\omega)$ be a compact symplectic $2n$-manifold. We assume that $(M,\omega)$
is quantizable; i. e. $[\omega]\in H^2(M,{\Bbb Z})$. Let $(L,D)$ be a prequantum
bundle \cite{nW92} over $M$.  We denote by $\{\psi_t\}_{t\in [0,1]}$  a Hamiltonian isotopy in $M$, 
with $\psi_0=\text{id}$. This family determines the set $\{X_t\}$ of vector fields by 
$$\frac{d\psi_t}{dt}=X_t\circ\psi_t.$$
By $f_t$ is denoted the corresponding normalized time-dependent Hamiltonian; that is,
$f_t$ is the function defined on $M$, such that $\omega(X_t,.)=-df_t$ and $\int_M f_t\omega^n=0$.

Denoting by $\Gamma(L)$ the space of $C^{\infty}$ sections of $L$, for each $t$ one defines the operator ${\cal P}_t\in\text{End}(\Gamma(L))$ by
 $${\cal P}_t(\sigma)=-D_{X_t}\sigma-2\pi if_t\sigma.$$
The differential equation
$$\frac{d\sigma_t}{dt}= {\cal P}_t(\sigma_t),\;\; \sigma_0=\sigma,$$
determines a family $\sigma_t$ of sections of $L$. In \cite[Corollary 5]{aV02} 
we have proved the 
following property: If
 $\psi_1=\text{id}$; i.e. $\psi$ is a loop in $\text{Ham}(M,\omega)$, then 
$\sigma_1=\kappa(\psi)\sigma$, where the constant $k(\psi)$ 
is given by the symplectic action around the curve $\{\psi_t(x) \}_t$
$$\kappa(\psi)=\text{exp}\Big(2\pi i\int_S\omega -2\pi i\int_0^1f_t(\psi_t(x))dt    \Big),$$
$x$ being any point of $M$ and $S$ being any $2$-chain in $M$ whose boundary is the nullhomologous
curve $\{\psi_t(x)\}_t$. Moreover $\kappa(\psi)$ depends only on the homotopy class $[\psi]$
of the loop $\psi$. That is, $\kappa$ is a representation 
$\kappa:\pi_1(\text{Ham}(M,\omega))\rightarrow U(1)$ \cite[Proposition 7]{aV02}.

\smallskip

Let $G$ be a compact connected Lie group. We denote by ${\cal O}:={\cal O}_{\eta}$ the
coadjoint orbit of $\eta\in{\frak g}^*$. ${\cal O}$ can be identified with $G/G_{\eta}$,
where $G_{\eta}$ is the stabilizer of $\eta$ for the coadjoint action. Given $A\in{\frak g}$,
by $X_A$ is denoted the vector field on ${\cal O}$ generated by $A$. The manifold ${\cal O}$ is equipped with
the symplectic structure $\omega$ defined by
$\omega(X_A(\nu),X_B(\nu))=\nu([A,B])\,$ \cite{aK76}. The map $\nu\in{\cal O}\mapsto -\nu\in{\frak g}^*$
is a moment map for the action of $G$ on ${\cal O}$; that is, $\iota_{X_A}\omega=-df_A$,
with $f_A\in C^{\infty}({\cal O})$ given by $f_A(\nu)=-\nu(A)$. If $A$ is a vector of 
${\frak Z}$, the center of ${\frak g}$, then $f_A$ is constant: $f_A(\nu)=\eta(A)$. So $f_A$ is
not normalized unless $\eta(A)=0$. Henceforth we assume that ${\frak Z}=0$.  

\smallskip

If the linear functional 
\begin{equation}\label{lambfunc}
\lambda:C\in{\frak g}_{\eta}=\{A\in{\frak g}\,|\,\eta([A,\,.])\equiv 0  \}\mapsto 2\pi i\eta(C)\in i{\Bbb R}
\end{equation}
is {\it integral}; that is, if there is a character $\Lambda:G_{\eta}\rightarrow U(1)$ whose
derivative is $\lambda$, then the orbit ${\cal O}$ is quantizable see \cite{bK70} \cite{aV02}.
A prequantum bundle $L$ over ${\cal O}=G/G_{\eta}$ is defined by
$L=G\times_{\Lambda}{\Bbb C}=(G\times {\Bbb C})/\simeq$, with
$(g,z)\simeq (gb^{-1}, \Lambda(b)z)$, for $b\in G_{\eta}$.

In this case each section $\sigma$ of $L$ determines a $\Lambda$-equivariant function $s:G\rightarrow{\Bbb C}$
by the relation
\begin{equation}\label{sfunction}
\sigma(gG_{\eta})=[g,s(g)].
\end{equation}
And given $A\in{\frak g}$ we denote by ${\cal P}_A$ the operator 
$$-D_{X_A}-2\pi if_A\in\text{End}(\Gamma(L)).$$
In \cite{aV02} is proved that the $\Lambda$-equivariant function associated
to ${\cal P}_A(\sigma)$ is $-R_A(s)$, where $R_A$ is the right invariant vector field on $G$ determined by $A$ and
 $s$ is the $\Lambda$-equivariant function associated to $\sigma$.

 \smallskip

If $\{g_t\}_{t\in[0,1]}$ is a smooth curve in $G$ with $g_0=e$, then
\begin{equation}\label{defpssit}
\{\psi_t:g'G_{\eta}\in G/G_{\eta} \mapsto g_tg'G_{\eta}\in G/G_{\eta}\},
\end{equation}
is a Hamiltonian isotopy.  This isotopy is generated by the vector fields $X_{A_t}$,
with $A_t=\Dot g_tg_t^{-1}$ (see \cite{aV02}).
The stabilizer $G_{\eta}$ contains a maximal torus of $G$ (see \cite{vG84}). Since the center
$Z(G)$ of $G$ is the intersection of all maximal tori of $G$, $\{\psi_t\}$ defines a loop
in $\text{Ham}({\cal O})$ if $g_1\in Z(G)$. From now on  
 we assume that $g_1\in Z(G)$.

One can consider the differential equation
$$\frac{d\sigma_t}{dt}={\cal P}_{A_t}\sigma_t,\;\; \sigma_t=\sigma.$$
If $s_t$ is the equivariant function associated to $\sigma_t$, by the above remark
$s_t$ satisfies $\Dot s(g_t)=-R_{A_t}(g_t)(s_t)$. But $\Dot g_t=R_{A_t}(g_t)\in T_{g_t}(G)$.
So $\Dot s_t(g_t)+\Dot g_t(s_t)=0$; that is,  the function $h:[0,1]\rightarrow {\Bbb C}$ defined
by $h(t)=s_t(g_t)$ is constant. Hence $s_1(g_1)=s_0(e)$.
$$\sigma_1(eG_{\eta})=[eG_{\eta}, s_1(e)]=[eG_{\eta},\Lambda(g_1)s_1(g_1)]=
\Lambda(g_1)[eG_{\eta}, s_0(e)]=\Lambda(g_1)\sigma_0(eG_{\eta}),$$
in other words, $\Lambda(g_1)=\kappa([\psi])$. We have the following Theorem
\begin{Thm}\label{carac}
Let $\{\psi_t\}_{t\in[0,1]}$ be the closed isotopy on ${\cal O}_{\eta}$ defined
 by $\psi_t(x)=g_t\cdot x$, where $g_t\in G$, \, $g_0=e$ and $g_1\in Z(G)$, if the  functional $2\pi i\eta$ is integral,
then
$\kappa([\psi])=\Lambda (g_1)$, where $\Lambda$ is the character of $G_{\eta}$
whose derivative is $2\pi i\eta$.
\end{Thm}

If $G_{\eta}$ is a maximal torus of $G$ there is another description of the 
action $\kappa([\psi])$ based on the Borel-Weil theorem \cite{D-K}. Here we quote the 
result of \cite{aV02}.
\begin{Prop} If $2\pi i\eta$ an integral character of ${\frak g}_{\eta}$ and $G_{\eta}$
is a maximal torus of $G$, then  the symplectic action around the loop (\ref{defpssit})
is 
\begin{equation}\label{kappaBorel}
\kappa([\psi])=\frac{\chi(\pi^*)(g_1)}{\text{dim}\,\pi},
\end{equation}
where $\pi$ is the irreducible representation of $G$ whose highest weight is $-2\pi i\eta$.
\end{Prop} 

%%%%%%%%%%%%%%%%%%%%%%%%%%%%%%%%%%%%%%%%%%%%%%%%%%%%%%%%%%%%%%%%%%%%%%%%%%%%%%%%
%%%%%%%%%%%%%%%%%%%%%%%%%%%%%%%%%%%%%%%%%%%%%%%%%%%%%%%%%%%%%%%%%%%%%%%%%%

\section{Symplectic action in flag manifolds}\label{SymFlag}

Let $G$ be a compact connected semisimple Lie group. As we said $G_{\eta}$
contains a maximal torus $T$. 
 Let us consider the decomposition of ${\frak g}_{\Bbb C}$
as direct sum of root spaces
$${\frak g}_{\Bbb C}={\frak t}_{\Bbb C}\oplus\bigoplus_{\alpha\in{\cal R}}{\frak g}_{\alpha}.$$
As $T\subset G_{\eta}$, $\eta$ vanishes over ${\frak g}_{\alpha}$. If
 $\check\alpha\in [{\frak g}_{\alpha},\,{\frak g}_{-\alpha}]$ is the coroot to $\alpha$, and
${\frak s}_{\alpha}:={\frak g}_{\alpha}\oplus {\frak g}_{-\alpha}\oplus [{\frak g}_{\alpha},\,{\frak g}_{-\alpha}]$,
it turns out that the complexification $({\frak g}_{\eta})_{\Bbb C}$ of ${\frak g}_{\eta}$
is generated by ${\frak t}_{\Bbb C}$ and the ${\frak s}_{\alpha}$'s for which $\eta(\check \alpha)=0$.
If we define 
$${\frak p}:={\frak t}_{\Bbb C}\oplus\bigoplus_{\eta(\check\alpha)\geq 0 }{\frak g}_{\alpha},$$
${\frak p}$ is a subalgebra which generates a parabolic subgroup $P$ of the
complexification $G_{\Bbb C}$ of $G$ \cite{F-H}, furthermore
$G_{\Bbb C}/P$ and $G/G_{\eta}$ can be identified as differential manifolds.
The element $\eta$ is said to be regular if $G_{\eta}=T$; that is, $\eta(\check\alpha)\ne 0$
for all $\alpha\in{\cal R}$.
In this case $G/G_{\eta}$ can be identified with $G_{\Bbb C}/B$, where $B$ is a Borel subgroup 
of   $G_{\Bbb C}$. Therefore the differential structure of 
${\cal O}_{\eta}$ is fixed by the set of roots $\alpha\in{\cal R}$ such that $\eta(\check\alpha)\ne 0$.
But in the symplectic structure of ${\cal O}_{\eta}$ are also involved the values of $\eta$.

If the functional $\lambda$ considered in   (\ref{lambfunc}) is integral (so
${\cal O}_{\eta}$ is quantizable)
 then $2\pi i\eta:{\frak t}\rightarrow i{\Bbb R}$
is the derivative of $\Lambda_{|T}$, where $\Lambda$ is a character of $G_{\eta}$. 
In this case, if $T\simeq{\Bbb T}^k$, there are integers $m_j$, $\,j=1,\dots, k$ such that
$\Lambda(A_1,\dots,A_k)=\prod A_j^{m_j}$, for $(A_1,\dots,A_k)\in{\Bbb T}^k$.
 Given  a smooth curve $\{g_t\in G\}_{t\in [0,1]}$ 
the symplectic action around the loop $\psi$ in $\text{Ham}({\cal O})$ 
defined in (\ref{defpssit}) is
$\kappa([\psi])=\Lambda(g_1)= \prod B_j^{m_j}$, assumed that $g_1\in Z(G)$ as
element of ${\Bbb T}^k$ is $g_1=(B_1,\dots,B_k) $.
So we have the following Theorem
\begin{Thm}\label{cardi}
If $\eta\in{\frak g}^*$ is integral, and the character $\Lambda$
of $G_{\eta}$ whose derivative is $2\pi i\eta$ is injective on $Z(G)$, then
$$\sharp(\big(\pi_1(\text{Ham}({\cal O}_{\eta}))\big)\geq \sharp(Z(G)).$$
\end{Thm}

\smallskip

{\it Note.} For a general coadjoint orbit ${\cal O}$ of $G$ McDuff and Tolman have
 proved the following property: If $G$ acts effectively on ${\cal O}$, the inclusion
$G\to\text{Ham}({\cal O})$ induces an injection on $\pi_1$ \cite{M-T}.
  
\smallskip

{\sc Example}

Given 
\begin{equation}\label{bfd}
{\bf d}=\big( \overbrace{p_1,\dots, p_1}^{n_1} ,\dots, \overbrace{p_k,\dots, p_k}^{n_k}  \big)\in{\Bbb R}^{n-1},
\end{equation}
with $0<n_1\leq n_2\leq\dots\leq n_k$,
we define $\eta\in{\frak su}(n)^*$ by 
$$\eta(Y)=i\sum_{j=1}^{n-1} d_jY_{jj},$$
where ${\bf d}=(d_1,\dots,d_{n-1})$. The stabilizer of $\eta$ for the 
coadjoint action of $G=SU(n)$
is $G_{\eta}=U(n_1)\times\dots\times U(n_k)\subset SU(n)$. The orbit 
$G/G_{\eta}$ is the flag manifold ${\cal F}_{\frak q}$ in ${\Bbb C}^n$, where 
${\frak q}$ is the partition $(n_1,\dots,n_k)$ of $n-1$. So the differential structure
of ${\cal O}_{\eta}$ is determined by the partition ${\frak q}$.

For 
$$Y=(B_1,\dots,B_k)\in\bigoplus_{j=1}^{k} {\frak u}(n_j)={\frak g}_{\eta}$$
 one has
$\eta(Y)=i\sum_{j=1}^kp_j\text{tr}(B_j).$ 
If $-2\pi p_j=:m_j\in{\Bbb Z}$ for $j=1,\dots, k$, then the
character $\Lambda:\prod U(n_j)\rightarrow U(1)$ defined by
$$\Lambda(A_1,\dots,A_k)=\prod_{j=1}^k(\text{det}(A_j))^{m_j}$$
has as derivative $2\pi i\eta$. That is, the coadjoint orbit ${\cal O}_{\eta}$ endowed with the natural symplectic
structure  is a quantizable manifold. Hence, given ${\frak q}=(n_1,\dots,n_k)$ a partition of $n-1$
and $\vec{m}\in{\Bbb Z}^k$, the pair $({\frak q},\,\vec{m})=:\vec d$ determines a quantizable symplectic
  manifold ${\cal O}_{\vec d}$, which is diffeomorphic to ${\cal F}_{\frak q}$

Given $z\in{\Bbb C}$ such that $z^n=1$, a curve $g_t\in SU(n)$ with $g_1=zI_n$
 defines
a closed Hamiltonian isotopy $_z\psi$ on ${\cal O}_{\vec d}$, by
  (\ref{defpssit}). 
   The corresponding symplectic action around the loop $_z\psi$ in $\text{Ham}({\cal O}_{\vec d})$
can be calculated by Theorem \ref{carac}
\begin{equation}\label{kappapsi0}
\kappa([_z\psi])=\prod_1^k z^{m_jn_j}=z^{{\frak q}\cdot\vec m}. 
\end{equation}
If ${\frak q}\cdot\vec m:=\sum_j n_jm_j$ and $n$ are relatively prime, then     
$\kappa([_z\psi])\ne\kappa([_y\psi])$, for $z\ne y$.
We have proved the following Proposition
\begin{Prop}\label{sun}
Let ${\cal O}_{\vec d}$ the coadjoint orbit of $SU(n)$ determined by 
$\vec d:=({\frak q},\vec m)$.
If ${\frak q}\cdot\vec m$ and $ n$ are relatively prime, then 
$\pi_1(\text{Ham}({\cal O}_{\vec d}))$ has at least $n$ elements.  
\end{Prop}

When $k=1$ in (\ref{bfd}); that is, $n_1=n-1$ and $m_1=:m$ then $G/G_{\eta}=SU(n)/U(n-1)={\Bbb C}P^{n-1}$. The
symplectic manifold ${\cal O}_{\eta}$ is ${\Bbb C}P^{n-1}$ equipped with a 
multiple of the Fubini-Study symplectic structure. If $g_1=zI_n$, with $z^n=1$
then by (\ref{kappapsi0})
\begin{equation}\label{kappapsi}
\kappa([_z\psi])=z^{-m}.
\end{equation}

In particular, when $m=-1$ the corresponding coadjoint orbit is ${\Bbb C}P^{n-1}$ endowed 
with the Fubini-Study symplectic structure $\omega_{FS}$. Therefore
\begin{Prop}\label{cpm-1wFS}
$\pi_1(\text{Ham}({\Bbb C}P^{n-1},\omega_{FS}))$ has at least $n$ elements.
\end{Prop}

{\it Remark 1.} It is known that 
$\text{Ham}({\Bbb C}P^{2},\omega_{FS})$ has the homotopy type of $PU(3)$ \cite{Polt}.
So $\#\pi_1(\text{Ham}({\Bbb C}P^{2},\omega_{FS}))=3$, and this result is consistent with
Proposition \ref{cpm-1wFS}.

\smallskip

 {\it Remark 2.} If $n=2$, then  ${\cal F}_{\frak q}={\Bbb C}P^1$, and
the value for $\kappa([_z\psi])$ obtained in (\ref{kappapsi})
can be also calculated using Proposition \ref{kappaBorel}.
Now  $-2\pi i\eta$
is the weight $(-m_1,0)$. Assumed that $ m<0$, if $\pi$ is the 
irreducible representation of $SU(2)$ whose highest weight is $-2\pi i\eta$, then
the  Weyl character formula \cite{G-W} gives
$$\chi_{\pi}(t_1,t_2)=\frac{\sum_{s\in{\cal S}_2}
\text{sig}\,s \;e^{s(-2\pi i\eta+\rho)}(t_1,\,t_2)}
{ \sum_{s\in{\cal S}_2}
\text{sig}\,s \;e^{s(\rho)}  (t_1,\,t_2) },$$
for $(t_1,t_2)\in U(1)\subset SU(2)$, $\,\rho$ being the weight $(1,0)$ and
${\cal S}_2$ the symmetric group of $2$ elements.
So
$$\chi_{\pi}(t_1,t_2)=\frac{t_1^{-m+1}- 
t_2^{-m+1} }{t_1-t_2    }=\sum_{k=0}^{-m}t_1^kt_2^{-m-k}$$
When $t_1=t_2=z$, 
$\chi_{\pi}(t_1,t_2)=(-m+1)z^{-m}$. On the other hand, Weyl dimension formula
gives $\text{dim}\,\pi=-m+1$.
 So, by Proposition \ref{kappaBorel}
$$\kappa([_z\psi])=\chi_{\pi}(z^{-1}I_2)/(-m+1)=z^{m}.$$
As $z^2=1$ this value agrees with   (\ref{kappapsi}).

\end{document}